\documentstyle{amsppt}

\pageheight{7.00in}
\NoBlackBoxes
\TagsOnRight

\define\bW {\partial\Omega}
\define\ep {\varepsilon}
\define\Ombar {\overline{\Omega}}
\define\Lip {\operatorname{Lip}}
\define\Hess {\operatorname{Hess}}
\define\wtN {\widetilde{\nabla}}

\define\RR {\Bbb R}
\define\pa {\partial}
\define\loc {\operatorname{loc}}

\magnification=1200

\document

\centerline{\bf LIPSCHITZ DOMAINS, DOMAINS WITH CORNERS,}
\vskip 10pt
\centerline{\bf AND THE HODGE LAPLACIAN
\footnote{2000 Math Subject Classification: Primary 31C12, 35B65, 35J25. 
Secondary 58J32, 42B20
\newline {\it Key words.} Hodge Laplacian, regularity, convexity, 
Lipschitz domains, manifolds with corners.
\newline M.\,M. was supported in part by NSF grant DMS-0400639
and the UMC Office of Research
\newline M.\,T. was supported in part by NSF grant DMS-0139726
\newline A.\,V. was partially supported by NSF grant DMS-0201092, a Clay 
Research Fellowship and a Fellowship from the Alfred P.\,Sloan Foundation}}

\vskip 15pt
\centerline{\smc Marius Mitrea, Michael Taylor, and Andr{\'a}s Vasy}

$$\text{}$$
{\smc Abstract}.  We define self-adjoint extensions of the Hodge Laplacian 
on Lipschitz domains in Riemannian manifolds, corresponding to either the 
absolute or the relative boundary condition, and examine regularity 
properties of these operators' domains and form domains.  
We obtain results valid for general Lipschitz domains, and stronger results 
for a special class of ``almost convex'' domains, which apply to domains
with corners.

$$\text{}$$
{\bf 1. Introduction}
\newline {}\newline

Let $\Omega$ be an open Lipschitz domain in a smooth, compact Riemannian 
manifold $M$, equipped with a metric tensor $g$, which we will assume is 
of class $C^2$. As is customary, let $d$, $\delta$ stand, respectively, for 
the operator of exterior differentiation and its adjoint. 
We use the Friedrichs method to define a self-adjoint extension
of the Hodge Laplacian $\Delta=-(d\delta+\delta d)$, with the absolute 
boundary condition (respectively, the relative boundary condition) on 
differential forms on $\Omega$, which we denote $-H=-H_A$ or $-H_R$.  We want 
to establish regularity properties of its domain $\Cal{D}(H)$ and of its form 
domain (which coincides with $\Cal{D}(H^{1/2})$).  We obtain a circle of 
results valid for general Lipschitz domains, and then some stronger results
valid for certain special classes of Lipschitz domains, including domains 
with corners.  These results extend some of the work in [MMT] and [M2].

To set things up, we define
$$
\aligned
X_A(\Omega)&=\{u\in L^2(\Omega,\Lambda^*):du,\delta u\in L^2(\Omega,
\Lambda^*),\nu\vee u|_{\bW}=0\}, \\
X_R(\Omega)&=\{u\in L^2(\Omega,\Lambda^*):du,\delta u\in L^2(\Omega,
\Lambda^*),\nu\wedge u|_{\bW}=0\}.
\endaligned
\tag{1.1}
$$
Here $\Lambda^*=\oplus_k \Lambda^k$, with $\Lambda^k$ denoting the $k$-th 
exterior exterior power of the tangent bundle of $M$, and $\nu$ is the unit 
conormal to $\bW$, a Lipschitz section of $\Lambda^1M|_{\bW}$; $\nu\vee u$ 
is an interior product and $\nu\wedge u$ an exterior product.  An important
ingredient in the proof that $X_A(\Omega)$ and $X_R(\Omega)$ are well
defined is the following result:
$$
\aligned
u,\delta u\in L^2(\Omega,\Lambda^*)&\Longrightarrow \nu\vee u\in 
H^{-1/2}(\bW,\Lambda^*), \\
u,du\in L^2(\Omega,\Lambda^*)&\Longrightarrow \nu\wedge u\in
H^{-1/2}(\bW,\Lambda^*),
\endaligned
\tag{1.2}
$$
established in (11.9) of [MMT]. Hereafter, $H^s$ will stand for 
the $L^2$-based Sobolev space of smoothness $s\in\Bbb{R}$, considered
either on $\Omega$ or on $\partial\Omega$. 
Also, $H^s(\Omega,\Lambda^*):=H^s(\Omega)\otimes\Lambda^*$, and so on, 
although in the sequel we shall occasionally drop the dependence of this, 
and other spaces, on the vector bundle.
 
If $(\cdot\,,\,\cdot)_{L^2(\Omega)}$ stands for the natural $L^2$-inner
product of forms in $\Omega$ (again, the subscript may be dropped in 
subsequent occurrences) then $X_A(\Omega)$ and $X_R(\Omega)$ are Hilbert 
spaces, with the inner product
$$
(du,dv)_{L^2}+(\delta u,\delta v)_{L^2}+(u,v)_{L^2}=Q(u,v)+(u,v)_{L^2},
\tag{1.3}
$$
where the last equality defines the sesqui-linear form $Q$. 
The Friedrichs extension method then yields self-adjoint operators $H_A$
and $H_R=-\Delta$ on $L^2(\Omega,\Lambda^*)$, with
$$
\aligned
\Cal{D}(H_A)=\{u\in X_A(\Omega):\,
X_A(\Omega)\ni v\mapsto Q(u,v)\text{ is }L^2\text{-bounded}\},
\\ (H_Au,v)_{L^2}=Q(u,v),\qquad\qquad\qquad\qquad
\endaligned
\tag{1.4}
$$
and $\Cal{D}(H_R)\subset X_R(\Omega)$ similarly defined.  As part
of the standard theory, one has
$$
\Cal{D}(H_A^{1/2})=X_A(\Omega),\quad 
\Cal{D}(H_R^{1/2})=X_R(\Omega).
\tag{1.5}
$$

In some cases, $X_A(\Omega)$ coincides with
$$
H^1_A(\Omega,\Lambda^*)=\{u\in H^1(\Omega,\Lambda^*):\nu\vee u|_{\bW}=0\},
\tag{1.6}
$$
with a similar result for $X_R(\Omega)$. This holds when $\bW$ is of 
class $C^2$, by a classical result of M.~Gaffney [G] and K.~Friedrichs [F].  
Also, if we write
$$
X_A(\Omega)=\bigoplus\limits_k X^k_A(\Omega),\quad
X_R(\Omega)=\bigoplus\limits_k X^k_R(\Omega),
\tag{1.7}
$$
where $X^k_A(\Omega)$ (respectively, $X^k_R(\Omega)$) consists of $k$-forms
in $X_A(\Omega)$ (respectively, in $X_R(\Omega)$), then standard regularity 
results for the Dirichlet and Neumann problems yield
$$
X^k_A(\Omega)=H^1_A(\Omega,\Lambda^k),\quad
X^k_R(\Omega)=H^1_R(\Omega,\Lambda^k),
\tag{1.8}
$$
for $k=0$ and $k=n$, where $n=\text{dim}\, \Omega$.  However, for general
Lipschitz $\Omega$ and $k\in [1,n-1]$, this identity fails.

$\text{}$ \newline
{\smc Example}.  Take $k=1,\ u=df$.  Then $du=0,\ \delta u=-\Delta f$, and
$\nu\vee u=\pa_\nu f$, so
$$
X^1_A(\Omega)\supset \{df:\,f\in H^1(\Omega),\,\Delta f\in L^2(\Omega),
\pa_\nu f=0\}.
$$
Now the regularity result
$$
f\in H^1(\Omega),\,\Delta f\in L^2(\Omega),\,\pa_\nu f=0\Longrightarrow 
f\in H^2(\Omega)
$$
is true if $\Omega$ is convex, or more generally satisfies a strong exterior
ball condition, but it fails for general Lipschitz $\Omega$.
\newline $\text{}$

It was shown in [MMT] that, for general Lipschitz $\Omega$,
$$
X_R(\Omega),\ X_A(\Omega)\subset H^{1/2}(\Omega,\Lambda^*).
\tag{1.9}
$$
A closer study of the example above shows that the exponent $1/2$ cannot
be improved in general. Furthermore, for a non-Lipschitz domain $\Omega$, 
elements in $X_A(\Omega)$, $X_R(\Omega)$ may fail to exhibit this 
critical amount of regularity. An example of a domain between two cones
with the same vertex and axis (thus not locally simply connected) is 
discussed in [CD2]. 

Various conditions ensuring the validity of (1.8) were given in [M2].  
These include a ``convexity'' hypothesis on $\Omega\subset M$, and a strong 
exterior ball hypothesis, in case $\Omega\subset
\RR^n$.  One of our main goals here is to extend that analysis, to include
a broader class of Lipschitz domains for which (1.8) is valid.  We define
a class of ``almost convex'' Lipschitz domains $\Omega$ in a compact
Riemannian manifold $M$.  We show that this class contains the class of 
Lipschitz domains $\Omega\subset M$ satisfying a uniform exterior ball 
condition, which in turn contains the class of compact manifolds with 
corners.  Furthermore we show that (1.8) holds for such almost convex 
domains.  Part of the interest in obtaining such a result is the potential
to extend the analysis of propagation of singularities in [Va] to the
setting of the wave equation $(\pa_t^2-\Delta)u=0$ when $\Delta$ is the 
Hodge Laplacian and $u=u(t,x)$ a differential form, on a manifold with 
corners, satisfying the absolute or relative boundary condition.
The regularity result (1.8) is also important in the variational
treatment of the Maxwell system in the class of forms of finite $L^2$-energy.
Cf., e.g., [CD2], [MM1], [MM2] for a discussion and references. 

The rest of the paper is structured as follows.  In \S{2} we present 
results on $\Cal{D}(H^{1/2})$ and on $\Cal{D}(H)$ valid for general Lipschitz
domains.  Some of these results are from [MMT], [MM1], and [MM2], and are
collected here for convenience.  Other results are new.  In \S{3} we 
introduce the notion of almost convexity, and show that it holds whenever 
the uniform exterior ball condition holds, and in particular that
domains with corners are almost convex.  In \S{4} we show that (1.8)
holds for almost convex domains, and establish further results on $\Cal{D}
(H)$ in this case.
\newline {}\newline

\noindent
{\it Acknowledgments.} The third author is grateful to Richard Melrose
for helpful discussions, and in particular for bringing the issue
of propagation of singularities for differential forms to his attention.

$$\text{}$$
{\bf 2. The Hodge Laplacian on Lipschitz domains}
\newline {}\newline

In this section we give further results on $X_A(\Omega)$ and on $\Cal{D}
(H_A)$ valid for general compact Lipschitz domains.  Note that since the Hodge 
star operator is its own inverse, up to sign, satisfies $\ast\Delta=\Delta\ast$
and has the mapping properties
$$
*:X_A(\Omega)\longrightarrow X_R(\Omega),\quad 
*:\Cal{D}(H_A)\longrightarrow \Cal{D}(H_R),
\tag{2.1}
$$
(plus a similar set with the roles of the subscripts $A$, $R$ reversed),  
it would suffice to investigate the absolute boundary condition.  
We also consider some other spaces:
$$
\aligned
X^k(\Omega)&=\{u\in L^2(\Omega,\Lambda^k):du,\delta u\in L^2(\Omega)\}, \\
X^k_b(\Omega)&=\{u\in X^k(\Omega):\nu\vee u,\nu\wedge u\in L^2(\bW)\}.
\endaligned
\tag{2.2}
$$
Recall that the result (1.2) makes $X^k_b(\Omega)$ well defined.  We have
the following trivial but occasionally useful observation:

\proclaim{Lemma 2.1} The spaces $X^k_A(\Omega), X^k_R(\Omega), X^k(\Omega)$,
and $X^k_b(\Omega)$ are all modules over $\Lip(\Ombar)$.
\endproclaim

The following result was established in Theorem 11.2 of [MMT].

\proclaim{Proposition 2.2} We have
$$
X^k_A(\Omega)\subset X^k_b(\Omega),
\tag{2.3}
$$
with an estimate
$$
\|u\|^2_{L^2(\bW)}\le C\bigl(\|du\|^2_{L^2(\Omega)}+\|\delta u
\|^2_{L^2(\Omega)}+\|u\|^2_{L^2(\Omega)}\bigr),\quad \forall\ 
u\in X^k_A(\Omega).
\tag{2.4}
$$
In fact, $X^k_A(\Omega)$ is a closed subspace of $X^k_b(\Omega)$.
\endproclaim

This result leads to the inclusion (1.9), when coupled with the following,
established in (11.20) of [MMT]:

\proclaim{Lemma 2.3} Given $u\in X^k_b(\Omega)$, we have, on $\Omega$,
$$
\aligned
u=\ &-d\Pi_{k-1}(\delta u)-\delta \Pi_{k+1}(du)-\Pi_k(Vu)
-\Cal{Q}_{k-1}(\delta u)-\Cal{R}_{k+1}(du) \\
&+\delta \Cal{S}_{k+1}(\nu\wedge u)-d\Cal{S}_{k-1}(\nu\vee u)
+R_{k+1}(\nu\wedge u)-R_{k-1}(\nu\vee u).
\endaligned
\tag{2.5}
$$
\endproclaim

Here $\Pi_k,\Cal{Q}_k$, and $\Cal{R}_k$ are integral operators on
forms on $\Omega$ (or even on $M$), and $\Cal{S}_k$ and $R_k$ are
layer potentials.  Also $V\in L^\infty(M)$.  Precise definitions of these
operators can be found in [MMT], particularly in (6.1)--(6.6) and (11.11).
We will state some of their mapping properties, for which we have further use 
below. We have
$$
\Pi_k:L^2(M)\longrightarrow H^2(M),\quad \Cal{Q}_k,\Cal{R}_k:
L^2(M)\longrightarrow H^1(M),
\tag{2.6}
$$
and
$$
\Cal{S}_k:L^2(\bW)\longrightarrow H^{3/2}(\Omega),\quad 
R_k:L^2(\bW)\longrightarrow H^{1/2}(\Omega).
\tag{2.7}
$$
Furthermore, for $p\in (1,\infty)$,
$$
\varphi\in L^p(\bW,\Lambda^*)\Longrightarrow \Cal{N}(d\Cal{S}_k\varphi),\,
\Cal{N}(\delta \Cal{S}_k \varphi),\,\Cal{N}(R_k \varphi)\in L^p(\bW),
\tag{2.8}
$$
where $\Cal{N}(\psi)$ denotes the nontangential maximal function associated
to a function or form $\psi$ on $\Omega$. At every boundary point $x\in\bW$, 
the latter is defined by $\Cal{N}(\psi)(x)=\text{sup}\,\{|\psi(y)|:\,
y\in\Omega,\,\text{dist}\,(x,y)<\kappa\,\text{dist}\,(y,\partial\Omega)\}$
for some fixed, sufficiently large $\kappa$. 

It should be mentioned that while (2.6) and the results on $R_k$ are
fairly straightforward, the results (2.7)--(2.8) on $\Cal{S}_k$ require 
the fundamental results of [Ca] and [CMM], and their extension to the setting
of potentials for variable coefficient operators, initiated in [MT] and 
carried out in the context needed here in Chapter 6 of [MMT].  It follows
from (2.5)--(2.7) that
$$
X^k_b(\Omega)\subset H^{1/2}(\Omega,\Lambda^k),
\tag{2.9}
$$
which together with (2.3) implies (1.9).  Furthermore, we have the following:

\proclaim{Corollary 2.4} Given $u\in X^k_b(\Omega)$, we have $u=T_1u+T_2u$,
with
$$
T_1u\in H^1(\Omega,\Lambda^k),\quad T_2u\in C^2_{\loc}(\Omega,\Lambda^k),
\quad \Cal{N}(T_2u)\in L^2(\bW).
\tag{2.10}
$$
\endproclaim
\demo{Proof}  Take $T_1u$ to be the sum of the first 5 terms on the right
side of (2.5), and take $T_2u$ to be the sum of the last 4 terms.  
The operators in (2.7) also map $L^2(\bW)$ to $C^2_{\loc}(\Omega)$.
\enddemo

We will improve (2.4) and also (2.10) (for $u\in X_A(\Omega)$)
later in this section, but for now we turn 
to other matters.  The following denseness result generalizes work in [CD1],
done there in the flat, three-dimensional Euclidean setting.

\proclaim{Proposition 2.5} For each $k\in\{0,\dots,n\}$, the space $C^2(
\Ombar,\Lambda^k)$ is dense in $X^k_b(\Omega)$.
\endproclaim
\demo{Proof} Via Lemma 2.1, we can assume $\Omega$ is a domain in $\RR^n$,
starlike about the origin (though with a variable coefficient, $C^2$ metric
tensor).  Then, given $u\in X^k_b(\Omega)$, write $u=T_1u+T_2u$ as in 
Corollary 2.4.  Certainly $T_1u$ is approximable in the $H^1$-norm, and
a fortiori in the $X^k_b$-norm, by elements of $C^2(\Ombar,\Lambda^k)$.
Furthermore, the dilates $w_r(x)=w(rx)$ for $r<1$ of $w=T_2u$ belong to
$C^2(\Ombar,\Lambda^k)$ and we have $w_r\rightarrow w, dw_r\rightarrow dw$,
and $\delta w_r\rightarrow \delta w$ in $L^2(\Omega)$ as $r\nearrow 1$,
and also $w_r|_{\bW}\rightarrow w|_{\bW}$ in $L^2(\bW)$; hence $w_r
\rightarrow w$ in $X^k_b(\Omega)$.
\enddemo

Proposition 2.5 is convenient for establishing some useful integration by
parts formulas. Throughout the paper, we let $dS$ denote the canonical surface
measure on $\partial\Omega$. Also, $\langle\cdot\,,\,\cdot\rangle$ 
stands for the {\it pointwise} inner product of forms. 

\proclaim{Proposition 2.6} Given $v\in X^k(\Omega)$ with $\delta v\in 
X^{k-1}_b(\Omega),\ dv\in X^{k+1}_b(\Omega)$, and $\varphi\in X^k_b(\Omega)$, 
we have
$$
(dv,d\varphi)+(\delta v,\delta\varphi)=-(\Delta v,\varphi)
+\int\limits_{\bW} \bigl[\langle \nu\vee dv,\varphi\rangle
-\langle \delta v,\nu\vee \varphi\rangle\bigr]\, dS.
\tag{2.11}
$$
In particular, if $v\in H^2(\Omega,\Lambda^k)$, $\varphi\in X^k_A(\Omega)$,
$$
(dv,d\varphi)+(\delta v,\delta\varphi)=-(\Delta v,\varphi)+
\int\limits_{\bW} \langle \nu\vee dv,\varphi\rangle \, dS.
\tag{2.12}
$$
\endproclaim
\demo{Proof} The identity (2.11) follows by adding up
$$
\aligned
(dv,d\varphi)&=(\delta dv,\varphi)+\int\limits_{\bW} \langle \nu\vee dv,
\varphi\rangle\, dS, \\
(\delta v,\delta\varphi)&=(d\delta v,\varphi)-\int\limits_{\bW}
\langle \delta v,\nu\vee\varphi\rangle\, dS.
\endaligned
$$
These identities, in turn, are easily justified by virtue of Proposition 2.5
and standard integration by parts formulas.
\enddemo

For applications below, it will be useful to complement Proposition 2.6 
with the following result.

\proclaim{Lemma 2.7} Given $\varphi\in X^k_b(\Omega)$ and $w\in C^2_{\loc}
(\Omega,\Lambda^k)$, satisfying
$$
\Delta w =f\in L^2(\Omega,\Lambda^k),\quad 
\Cal{N}(w),\ \Cal{N}(dw),\ \Cal{N}(\delta w)\in L^2(\bW),
\tag{2.13}
$$
with the boundary values taken in $L^2(\bW)$, we have
$$
(dw,d\varphi)+(\delta w,\delta\varphi)=-(f,\varphi)-\int\limits_{\bW}
\langle \nu\vee dw,\varphi\rangle\, dS+\int\limits_{\bW}
\langle \delta w,\nu\vee \varphi\rangle\, dS.
\tag{2.14}
$$
\endproclaim
\demo{Proof} Take a sequence $\Omega_\ell\subset\subset\Omega$ such that
$\Omega_\ell\nearrow\Omega$ in a nice fashion and denote by $\nu_\ell$, 
$dS$, respectively, the unit conormal and surface measure on 
$\partial\Omega_\ell$. Then $w|_{\Omega_\ell}\in
C^2(\Ombar,\Lambda^k)$ and $\varphi|_{\Omega_\ell}\in H^1(\Omega_\ell,
\Lambda^k)\subset X^k_b(\Omega)$, so (2.12) applies, to give
$$
\aligned
(dw,d\varphi)_{L^2(\Omega_\ell)}+(\delta w,\delta\varphi)_{L^2(\Omega_\ell)}
=-(f,\varphi)_{L^2(\Omega_\ell)}&-\int\limits_{\pa\Omega_\ell}
\langle \nu_\ell\vee dw,\varphi\rangle\, dS_\ell \\
&+\int\limits_{\pa\Omega_\ell} \langle \delta w,\nu_\ell\vee \varphi
\rangle\, dS_\ell.
\endaligned
\tag{2.15}
$$
It is elementary that the left side of (2.15) converges to the left side of
(2.14) as $\ell\rightarrow \infty$.  Now the boundary behaviors of $\varphi$
and $w$ given by hypothesis also yield convergence of the right side of 
(2.15) to the right side of (2.14) as $\ell\rightarrow\infty$, so we have
the lemma.
\enddemo

It is of interest to look at the Dirac-type operator
$$
D_A=d+\delta,\quad \Cal{D}(D_A)=X_A(\Omega),
\tag{2.16}
$$
and its counterpart $D_R=d+\delta,\ \Cal{D}(D_R)=X_R(\Omega)$.  The 
Friedrichs construction of $H_A$ and $H_R$ entails
$$
H_A=D_A^* D_A,\quad H_R=D_R^* D_R.
\tag{2.17}
$$
In light of this, it is valuable to have the following, which follows
from Proposition 6.1 and Theorem 6.2 of [MM2]:

\proclaim{Proposition 2.8} The operators $D_A$ and $D_R$ are self-adjoint.
\endproclaim

Hence we have
$$
H_A=D_A^2,\quad H_R=D_R^2,
\tag{2.18}
$$
and consequently
$$
\Cal{D}(H_A)=\{u\in X_A(\Omega):(d+\delta)u\in X_A(\Omega)\}.
\tag{2.19}
$$
Note that $H_A$ takes $k$-forms to $k$-forms, and we can write
$$
H_A=\bigoplus\limits_k H_{A,k},\quad
\Cal{D}(H_{A,k})=\Cal{D}(H_A)\cap L^2(\Omega,\Lambda^k).
\tag{2.20}
$$
We see that
$$
\Cal{D}(H_{A,k})=\{u\in X^k_A(\Omega):\,du\in X^{k+1}_A(\Omega),\,
\delta u\in X^{k-1}_A(\Omega)\}.
\tag{2.21}
$$
In particular,
$$
u\in \Cal{D}(H_{A,k})\Longrightarrow \nu\vee du=0\ \text{ on }\ \bW.
\tag{2.22}
$$
Membership of $u$ to $\Cal{D}(H_{A,k})$ also entails $\nu\vee \delta u=0$ 
on $\bW$, though this is automatic from $\nu\vee u=0$ and (6.15) of [MMT].  
Another consequence of (2.21) is that
$$
u\in \Cal{D}(H_{A,k})\Longrightarrow \delta du,\,d\delta u\in 
L^2(\Omega,\Lambda^k).
\tag{2.23}
$$

The following result gives important additional information on $\Cal{D}
(H_A)$.  Related work, in the more general context of $L^p$ with $p$ close 
to $2$, is given in \S{5} of [MM1]; cf. also [M1]. Let us also note here 
that, as far as the optimality of the range of possible $p$'s is concerned, 
the case of three-dimensional manifolds is best understood at the moment. 
In this context, sharp estimates on Sobolev-Besov spaces for the Hodge 
Laplacian on Lipschitz domains have been recently proved in [M3]. 

\proclaim{Proposition 2.9} There exists $p=p(\Omega)>2$ with the
following property.
Let $u\in\Cal{D}(H_{A,k})$.  Assume $1\le k\le 
n-1$ (since otherwise stronger results hold).  Then $u=v-w$ with
$$
\aligned
v\in H^2(\Omega),\quad w\in C^2_{\loc}(\Omega),\quad
(\Delta-1)w=0,\\
\Cal{N}(w),\,\,\Cal{N}(dw),\,\,\Cal{N}(\delta w)\in L^{p}(\bW).\qquad
\endaligned
\tag{2.24}
$$
Also, the boundary values of $w$, $dw$ and $\delta w$ exist in $L^p(\bW)$.
\endproclaim
\demo{Proof} Let $\Cal{O}$ be an open neighborhood of $\overline{\Omega}$.
Given $F\in L^2(\Omega,\Lambda^k)$, extend
$F$ to $\Cal{O}$ and solve for $v$:
$$
(\Delta-1) v=F,\quad v\in H^2_{\loc}(\Cal{O}).
\tag{2.25}
$$
Then
$$
f=\nu\vee v,\ g=\nu\vee dv\in L^p(\bW)
\tag{2.26}
$$
for some $p>2$ depending only on $n=\text{dim}\, \Omega$.
Furthermore, by (6.15) of [MMT], we have
$$
f\in L^{p,\delta}_{\text{tan}}(\bW,\Lambda^{k-1}TM),
\tag{2.27}
$$
a space defined by (5.2) of [MMT].  Hence, by (a simple variant of)
Theorem 5.1 of [MMT], there exists $w\in C^2_{\loc}(\Omega,\Lambda^k)$ 
such that (possibly with smaller $p>2$)
$$
\Cal{N}(w),\quad \Cal{N}(dw),\quad \Cal{N}(\delta w)\in L^{p}(\bW),
\tag{2.28}
$$
and
$$
(\Delta-1) w=0,\quad \nu\vee w=f,\quad \nu\vee dw=g.
\tag{2.29}
$$
Set $u=v-w$.  We see that
$$
du=dv-dw\in L^2(\Omega),\quad \delta u=\delta v-\delta w\in L^2(\Omega),
\tag{2.30}
$$
and that
$$
\nu\vee u=\nu\vee v-\nu\vee w=f-f=0,
\tag{2.31}
$$
so
$$
u\in X^k_A(\Omega).
\tag{2.32}
$$
Furthermore, given $\varphi\in X^k_A(\Omega)$, we have from (2.12) and
(2.14) that
$$
\align
(dv,d\varphi)+(\delta v,\delta\varphi)&=-(F+v,\varphi)-\int\limits_{\bW}
\langle g,\varphi\rangle\, dS, \tag{2.33} \\
(dw,d\varphi)+(\delta w,\delta\varphi)&=-(w,\varphi)-\int\limits_{\bW}
\langle g,\varphi\rangle\, dS, \tag{2.34}
\endalign
$$
so that subtracting (2.34) from (2.33) yields
$$
(du,d\varphi)+(\delta u,\delta\varphi)=-(F+u,\varphi),\quad \forall\ 
\varphi\in X^k_A(\Omega).
\tag{2.35}
$$
Thus $u\in\Cal{D}(H_{A,k})$ and $H_{A,k}u+u=-F$.  

Since this works for arbitrary $F\in L^2(\Omega,\Lambda^k)$ and 
since $1+H_{A,k}$ has a bounded inverse on $L^2(\Omega,\Lambda^k)$, 
this proves the proposition. 
\enddemo

We are now ready for the advertised improvements on (2.4) and (2.10).
First, recall the operators $T_1$, $T_2$ introduced in Corollary 2.4. 

\proclaim{Proposition 2.10} There exists $p=p(\Omega)>2$ such that
$$
u\in X_A(\Omega)\Longrightarrow u\bigr|_{\bW}\in L^p(\bW).
\tag{2.36}
$$
Furthermore, we have $u=T_1u+T_2u$ with $T_1u\in H^1(\Omega,\Lambda^*),\
T_2u\in C^2_{\loc}(\Omega,\Lambda^*)$, and
$$
\Cal{N}(T_2u)\in L^p(\bW).
\tag{2.37}
$$
\endproclaim
\demo{Proof} From Proposition 2.8 and (2.18) we see that any $u\in 
X_A(\Omega)$ can be written
$$
u=(d+\delta)v+w,\quad v\in\Cal{D}(H_A),\ \ w\in\text{Ker}\, H_A.
\tag{2.38}
$$
In fact, this is a manifestation of the fact that $L^2$-Hodge decompositions 
of forms are valid in any Lipschitz domain; cf. [MMT] and [MM1]. 
Thus (2.36) follows from (2.24).  Then (2.37) follows from another application
of (2.8), given the formula for $T_2u$.
\enddemo

Theorem 7.4 of [MT2] provides the following complement to (2.7):
$$
\Cal{S}_k:L^p(\bW)\longrightarrow B^{p,p^\#}_{1+1/p}(\Omega),\quad 
1<p<\infty,
\tag{2.39}
$$
where the target space is a Besov space and $p^\#:=\text{max}\,\{p,2\}$. 
One has a corresponding result for $R_k$ on $L^p(\bW)$.  
Hence (2.37) can be complemented by
$$
T_2 u\in B^{p,p}_{1/p}(\Omega),
\tag{2.40}
$$
for some $p>2$. Consequently, we have the following. 

\proclaim{Corollary 2.11}
There exists $p=p(\Omega)>2$ such that  
$$
X_A(\Omega),\,X_R(\Omega)\subset B^{p,p}_{1/p}(\Omega,\Lambda^*). 
\tag{2.41}
$$
\endproclaim

Let us remark that this regularity result is in the nature of best possible
in the class of Lipschitz domains. 
Indeed, for $\omega\in (0,\pi)$ we let $\Omega_\omega$ be the two-dimensional 
domain which coincide with the sector 
$\{z\in\Bbb{C}:\,|\text{arg}\,z|<\omega\}$ near the origin, 
$w(z)=\text{Re}\,(z^{\pi/2\omega})$, and finally set $u=dw$, suitably 
truncated near the origin. Then $u\in X_R(\Omega_\omega)$ and 
$$
u\in B^{p,p}_{1/p}(\Omega_\omega,\Lambda^1)\Longleftrightarrow 
p<\frac{2\omega}{2\omega-\pi}.
\tag{2.42}
$$
Note that $\omega\nearrow\pi$ forces $\frac{2\omega}{2\omega-\pi}\searrow 2$, 
justifying the claim about the sharpness of (2.41). This example can be 
modified to work in higher dimensions by adding extra dummy variables. 

Closer inspection of the situation described above reveals that, nonetheless, 
$u\in H^1(\Omega_\omega,\Lambda^1)$ whenever $0<\omega<\frac{\pi}{2}$, 
in which case $\Omega_\omega$ is geometrically {\it convex}. 
This type of phenomenon is examined in the greater detail in Sections 3-4.

$$\text{}$$
{\bf 3. Almost convex domains and domains with corners}
\newline {}\newline

To set up our first definitions, we assume we have a nested family of 
$C^2$ domains $\Omega_\ell\nearrow\Omega$.  
We take a neighborhood $U$ of $\bW$ and assume $\pa\Omega_\ell\subset U$
for all $\ell$.  We assume each $\Omega_\ell$ has a $C^2$ defining function 
$\rho_\ell$, defined on $U$, strictly negative on $\Omega_\ell\cap U$ and 
vanishing on $\pa\Omega_\ell$, satisfying
$$
C_1^{-1}\le \|d\rho_\ell(x)\|\le C_1,\quad \forall\ x\in \bW_\ell,
\tag{3.1}
$$
for some $C_1\in (0,\infty)$.  This ensures that each $\Omega_\ell$ is
Lipschitz, with Lipschitz constant independent of $\ell$.
The norm in (3.1) is defined by the metric tensor, but of course the 
condition (3.1) is independent of choice of metric tensor. The {\it Hessian}
of $\rho_\ell$ is defined by $\Hess(\rho_\ell)=\nabla d\rho_\ell$,
where $\nabla$ is the Levi-Civita connection of $g$.
In local coordinates $x$ this takes the form
$$
\Hess(\rho_\ell)=\sum\limits_{i,j} \left(\frac{\pa^2 \rho_\ell}{\pa x_i \pa x_j}
-\sum_k \Gamma^k_{ij} \frac{\pa \rho_\ell}{\pa x_k}\right)\,
dx_i\, dx_j,
\tag{3.2}
$$
where $\Gamma^k_{ij}$ are the Christoffel symbols of $g$.
Note that by (3.1), the second term is uniformly bounded in $\ell$
over compact subsets of the coordinate patch, while the
first term is the Hessian, $\Hess_x(\rho_\ell)$,
of $\rho_\ell$ with respect to the Euclidean metric on the coordinate patch.

Our hypothesis of almost convexity is:
$$
\Hess(\rho_\ell)\geq -C_2 g,
\tag{3.3}
$$
as quadratic forms on $T\bW_\ell$,
for some $C_2\in (0,\infty)$, independent of $\ell$.
In view of (3.2) and (3.1), an equivalent formulation is the following.
Cover $U$ by a finite number of coordinate systems $(O^i,x^i)$, let
${\Cal O}^i\subset O^i$ satisfy $\overline{{\Cal O}^i}\subset O^i$,
and $\cup_i {\Cal O}^i\supset U$, and assume that in
each of these local coordinates $x$, over ${\Cal O}^i$,
$$
\sum\limits_{i,j} \frac{\pa^2 \rho_\ell}{\pa x_i \pa x_j} \xi_i \xi_j
\ge -C_2 \sum\limits_i \xi_i^2,\quad \text{whenever }\ \rho_\ell=0
\ \text{and}\ 
\sum\limits_i \frac{\pa \rho_\ell}{\pa x_i} \xi_i=0,
\tag{3.4}
$$
for some (perhaps different) $C_2\in (0,\infty)$, independent of $\ell$. 

This notion is independent of the choice of metric tensor, since for two
Riemannian metrics $g$ and $g'$, by (3.2), $\nabla_g-\nabla_{g'}$ is
a zeroth order differential operator from $T^*M$ to $T^*M\otimes T^*M$,
so $\Hess_g \rho_\ell-\Hess_{g'}\rho_\ell$ is uniformly bounded (from above
as well as below) as a quadratic form by (3.1).

Alternatively, if we take the equivalent definition (3.4), which
is clearly metric, but not coordinate, independent, then we can see
directly that almost convexity is independent of the coordinate systems
chosen. In fact, take two coordinate systems
$x=(x_1,\dots,x_n)$ and $y=(y_1,\dots,y_n)$.  By the chain rule, we have
$$
\Hess_x(\rho_\ell)=\sum\limits_{i,j,r,s,m,k} \frac{\pa y_r}{\pa x_i}
\frac{\pa}{\pa y_r}\Bigl(\frac{\pa y_s}{\pa x_j}
\frac{\pa\rho_\ell}{\pa y_s}\Bigr)\, \frac{\pa x_i}{\pa y_k}
\frac{\pa x_j}{\pa y_m}\, dy_k\, dy_m=I_\ell+II_\ell,
\tag{3.5}
$$
where the terms $I$ and $II$ arise by the derivative $\pa/\pa y_r$ falling
on $\pa y_s/\pa x_j$ or $\pa\rho_\ell/\pa y_s$, respectively.
By (3.1), the coefficient of each term $dy_k\, dy_m$ in $I_\ell$ is 
uniformly bounded in $\ell$, so we have $-C_3 g\le I_\ell\le C_3g$
for a suitable constant $C_3$.  On the other hand,
$$
II_\ell=\sum\limits_{k,m} \frac{\pa^2 \rho_\ell}{\pa y_k\pa y_m}\, 
dy_k\, dy_m=\Hess_y(\rho_\ell),
\tag{3.6}
$$
so the coordinate independence of the definition of almost convexity
is established.  

We state another characterization of almost convexity. Let $l_\ell:T\bW_\ell
\otimes T\bW_\ell\to\RR$ denote the real-valued second fundamental form
of $\bW_\ell$. Recall that $l_\ell$ depends on the choice of a unit normal
vector field at $\bW_\ell$; we normalize it using the outward pointing
normal vector $n_\ell$ to $\bW_\ell$. Then $l_\ell$ is given by $l_\ell(u,v)
=g(\nabla_u n_\ell,v)$, $u,v\in T_p\bW_\ell$,
where $\nabla$ is the covariant derivative in $M$.
This can be rephrased
in terms of the outward pointing conormal, $\nu_\ell$.
Namely, using $g(n_\ell,v)=0=\nu_\ell(v)$ for $v\in T\bW_\ell$, we deduce that
$l_\ell(u,v)=(\nabla_u \nu_\ell)(v)$. Then
$\nu_\ell=\frac{d\rho_\ell}{\|d\rho_\ell\|}$, and
$l_\ell=\frac{\Hess\rho_\ell}{\|d\rho_\ell\|}$. Thus, almost convexity
is equivalent to requiring that $l_\ell$ be bounded below, uniformly in $\ell$.

We can compare the notion of almost convexity with that of 
$k$-convexity in the sense of [M2], which requires
$$
\sum\limits_{i,j} (\rho_\ell)_{ij}
\langle w_i\wedge u,w_j\wedge u\rangle\ge 0
\tag{3.7}
$$
on $\pa\Omega_\ell$,
where $\{w_i:1\le i\le n\}$ is a local orthonormal frame field for $T^*U$,
$u$ is an $k$-form ($1\le k\le n$) satisfying $\nu\wedge u=0$, and 
$\langle\cdot\,,\,\cdot\rangle$ is the (pointwise) inner product on 
$\Lambda^{k+1}T^*U$ induced by the metric tensor $g$.  
Also, $(\rho_\ell)_{ij}$ stands for the expression 
$\frac12\Bigl\{\frac{\partial^2\rho_\ell}{\partial w^i\partial w^j}
+\frac{\partial^2\rho_\ell}{\partial w^j\partial w^i}\Bigr\}$, 
where $\{\partial/\partial w^i\}_i$ form the (local) orthonormal basis of $TM$ 
dual to $\{w^i\}_i$. The proof of Proposition 3.2 of [M2] 
shows that if $\Omega$ is almost convex, then one has, in place of (3.7),
$$
\sum\limits_{i,j} (\rho_\ell)_{ij}
\langle w_i\wedge u,w_j\wedge u\rangle\ge -C\langle u,u\rangle,
\tag{3.8}
$$
for some $C\in (0,\infty)$, for such $u$ as above, and for all $k\le n$.
 
We now discuss some important special classes of almost convex domains.

We say that a Lipschitz domain $\Omega\subset M$ satisfies a local exterior 
ball condition, henceforth referred to as LEBC, 
if for every boundary point $x_0\in\partial\Omega$
there exists a coordinate patch ${\Cal O}$ which contains $x_0$ and
which satisfies the following two conditions. 

First, there exists a Lipschitz function $\varphi:{\RR}^{n-1}\to{\RR}$ with 
$\varphi(0)=0$ and such if $D$ is the domain above the graph of $\varphi$
then $D$ satisfies the standard Euclidean uniform exterior ball condition
(UEBC).  Second, it is assumed that there exists a diffeomorphism 
$\Upsilon$ mapping ${\Cal O}$ onto the unit ball $B_1(0)$ in ${\RR}^n$ 
and such that $\Upsilon(x_0)=0$, $\Upsilon({\Cal O}\cap\Omega)=B_1(0)\cap D$, 
$\Upsilon({\Cal O}\setminus{\Ombar})=B_1(0)\setminus\overline{D}$. 

\proclaim{Proposition 3.1} If the Lipschitz domain $\Omega\subset M$ 
satisfies a LEBC then it is almost convex. 
\endproclaim
\demo{Proof}We need to construct a nested sequence of approximating 
domains $\Omega_\ell$ of bounded Lipschitz character, 
which have $C^2$ defining functions $\rho_\ell$ 
with the properties described above.

The construction is local in nature and, given that the original domain
satisfies a LEBC, there is no loss of generality in assuming that 
a system of local coordinates $\{x_i\}_i$ has been selected for which 
$x_0=0$, $g$ is the Euclidean metric in this system of coordinates,  
and $\Omega$ is the domain above the graph of a Lipschitz function 
$\varphi:{\RR}^{n-1}\to{\RR}$ with $\varphi(0)=0$ satisfying
the following property. There exists $C_0>0$ such that for a.e.  
$a\in{\RR}^{n-1}$ and $\forall\,v\in 
{\RR}^{n-1}$, $|v|\leq C_0$, there holds 
$$
2\varphi(a)-\varphi(a+v)-\varphi(a-v)\leq C_0|v|^2. 
\tag{3.9} 
$$
That the latter condition can be assumed is a consequence of
our definition of LEBC and Lemma 6.3 in [M2]. 
In this scenario, following the construction in \S{6} of [M2], we 
take, for each $\ell\ge 1$, 
$$
\rho_\ell(x):=\varphi_\ell(x')-x_n,\qquad x=(x',x_n)\in{\RR}^n. 
\tag{3.10} 
$$
Above, for each $x'\in{\RR}^{n-1}$, 
$$
\varphi_\ell (x'):=\frac{C}{\ell} 
+\int\limits_{{\RR}^{n-1}}\Phi_\ell(y')\varphi(x'-y')\,dy',
\tag{3.11} 
$$
where $C>0$ is a fixed constant (to be specified below), 
$\Phi\in C^\infty({\RR}^{n-1})$, $0\leq\Phi\leq 1$, $\Phi\equiv 0$ 
for $|x'|>1$, $\int_{{\RR}^{n-1}}\Phi\,dx'=1$ and, as is customary,  
$\Phi_\ell(x'):=\ell^{n-1}\Phi(\ell x')$. If we now set 
$$
\Omega_\ell:=\{x:\,\rho_\ell(x)<0\}=\{(x',x_n):\,\varphi_\ell(x')<x_n\}, 
\tag{3.12} 
$$
it follows that $\partial\Omega_\ell\in C^\infty$.  Furthermore, if the 
constant $C$ in (3.11) is sufficiently large, then the family (3.12)
is nested, ${\Ombar}_\ell\subseteq\Omega$ and 
$\cup_{1<\ell<\infty}\Omega_\ell=\Omega$. More specifically, if 
$C>\|\nabla\varphi\|_{L^\infty}$ then 
$$
\frac{d}{d\mu}\Bigl[\varphi_{1/\mu}(x')\Bigr]=C-
\int\limits_{{\RR}^{n-1}}\Phi(z')\,z'\cdot(\nabla\varphi)(x'-\mu z')\,dz'>0
\tag{3.13}
$$
which ensures that mapping $\ell\mapsto\varphi_\ell(x')$ is decreasing. 
Thus, $\varphi_\ell(x')\searrow\varphi(x')$ as $\ell\nearrow\infty$. 
In addition, by virtue of (3.9) --\,a feature also inherited by 
each $\varphi_\ell$\,-- the domain $\Omega_\ell$ satisfies a UEBC with a 
constant independent of $\ell\in(1,\infty)$. Thanks to Lemma 6.3 in [M2], 
this last property further entails the existence of some $C>0$ such that 
$$
{\Hess}(\varphi_\ell)\geq -C
\tag{3.14} 
$$
uniformly in $\ell\in(1,\infty)$. Here ${\Hess}(\varphi_\ell)$ is the Hessian 
of $\varphi_\ell$ in the coordinates $\{x_i\}_{1\leq i\leq n-1}$, viewed as a
symmetric $(n-1)\times(n-1)$ matrix. 

It follows that for {\it all} vectors $t=(t',t_n)\in{\RR}^n$,   
$$  
{\Hess}_{x}(\rho_\ell)\,t\,\cdot\,t=
{\Hess}(\varphi_\ell)\,t'\,\cdot\, t'\geq -C|t'|^2\geq -C|t|^2, 
\tag{3.15} 
$$
uniformly in $\ell$. This is one of the two conditions
we set to check (recall that we have chosen $g$ to be the
Euclidean metric in the system of coordinates $\{x_i\}_i$)
The remaining one, (3.1), is easily seen from (3.11). Indeed,
$\|d\rho_\ell(x)\|\approx(1+|\nabla\varphi_\ell(x)|^2)^{1/2}$
and $|\nabla\varphi_\ell(x)|\leq |\nabla\varphi(x)|$, uniformly in $\ell$. 
\enddemo

$\text{}$ \newline
{\smc Remark}. What (3.9) says is that, in the approximation scheme 
$\Omega_\ell\nearrow\Omega$ we have constructed for a 
domain $\Omega$ satisfying a LEBC, the Hessians of the defining functions 
$\rho_\ell$ for $\Omega_\ell$ are bounded from below on the {\it entire} 
tangent space to $M$, uniformly in $\ell$, rather than just on 
$T\partial\Omega_\ell$ as required for almost convex domains. (In fact,
if $\rho_\ell$ satisfy the almost convexity hypotheses (3.1) and (3.3),
they can be
replaced by some $\tilde\rho_\ell=F_\ell(\rho_\ell)$ so that for
$\tilde\rho_\ell$, (3.1) still
holds, and the lower
bound (3.3) holds on the entire tangent space.)
This gives a heuristic explanation as to why domains with LEBC 
happen to be almost convex. 
\newline $\text{}$

A Lipschitz domain $\Omega\subset M$ is said to be a domain with corners
provided that each $p\in\bW$ has a neighborhood $\Cal{O}$ on which there
are coordinates $x_1,\dots,x_n$ such that $\Ombar\cap\Cal{O}$ is defined by
$$
x_j\ge 0\ \text{ for }\ 1\le j\le m,
\tag{3.16}
$$
for some $m\in\{1,\dots,n\}$.  The following is apparent.

\proclaim{Proposition 3.2} Every domain with corners satisfies a LEBC, and 
hence is almost convex.
\endproclaim

$\text{}$ \newline
{\smc Remark}. It is an interesting exercise to prove {\it directly} that
any domain with corners is almost convex. This can be shown locally, in 
particular in some local coordinate system $(x_1,...,x_n)$ where it may be
assumed that $\Omega$ is given by the system of inequalities (3.16) 
plus the requirement that all $x_j$'s are bounded. Then one can define
the family $\Omega_\ell$ as $\{\rho_\ell<0\}$ where 
$\rho_\ell=\tilde{\rho}_\ell/\|d\tilde{\rho}_\ell\|$ and 
$\tilde{\rho}_\ell=\ell-x_1\cdots x_m$. Thus, clearly, (3.1) holds. 
To verify (3.4) note that the condition 
$\sum\limits_i (\pa \rho_\ell/\pa x_i) \xi_i=0$ entails 
$\sum\limits_i \xi_i/x_i=0$. On the other hand, the semi-boundedness
condition on the Hessian amounts to checking that 
$$
\sum\limits_{i\not=j} \frac{-1}{\ell}\frac{\xi_i\xi_j}{x_ix_j}
\ge -C|\xi|^2\|d\tilde{\rho}_\ell\|
$$ 
whenever the vector 
$\sum\limits_i \xi_i(\partial/\partial x_i)$ is tangent to the
zero set of $\rho_\ell$. However, for such a vector, 
$$
0=\Bigl(\sum_i\frac{\xi_i}{x_i}\Bigr)\Bigl(\sum_j\frac{\xi_j}{x_j}\Bigr)
=\sum_{i,j}\frac{\xi_i\xi_j}{x_ix_j}
=\sum_{i}\frac{\xi_i^2}{x_i^2}+\sum_{i\not=j}\frac{\xi_i\xi_j}{x_ix_j}. 
\tag{3.16}
$$
As the first term in the rightmost expression is non-negative, the Hessian
condition follows (with $C=0$).

$$\text{}$$
{\bf 4. The Hodge Laplacian on almost convex domains}
\newline {}\newline

We aim to prove the following.

\proclaim{Theorem 4.1} If $\Omega$ is an almost convex domain, then
$$
\Cal{D}(H^{1/2}_A)=H^1_A(\Omega,\Lambda^*)=\bigl\{u\in H^1(\Omega,\Lambda^*):
\nu\vee u\bigr|_{\bW}=0\bigr\},
\tag{4.1}
$$
and
$$
\Cal{D}(H^{1/2}_R)=H^1_R(\Omega,\Lambda^*)=\bigl\{u\in H^1(\Omega,\Lambda^*):
\nu\wedge u\bigr|_{\bW}=0\bigr\}.
\tag{4.2}
$$
\endproclaim

From this and Proposition 3.2 we may therefore readily conclude the following. 

\proclaim{Corollary 4.2} The identities (4.1)-(4.2) hold for any 
domain $\Omega$ with corners.
\endproclaim

A quick sketch of the proof of Theorem 4.1 is as follows.  
For a given form $u$, we follow the proof of Theorem 5.1 in [M2], 
noting that the only place where a modification is needed is equation (5.9).  
There, the boundary term
$$
\sum\limits_{i,j} \int\limits_{\pa\Omega_\ell} (\rho_\ell)_{ij}
\langle w_i\wedge v_\ell,w_j\wedge v_\ell\rangle\, dS_\ell
$$
(in [M2], the index $\ell$ is actually denoted by $\mu$) 
can be dropped from a subsequent estimate since it is non-negative by the 
convexity assumption (3.7).  However, it suffices if one can estimate this
term from below by $-C\|v_\ell\|_{H^1(\Omega_\ell)}\|v_\ell
\|_{L^2(\Omega_\ell)}$ for some $C>0$, for then the effect on (5.11) in [M2]
is the same as the remainder term (5.10), which is already controlled.
But this follows from (3.8) combined with
$$
\|v_\ell\|^2_{L^2(\pa\Omega_\ell)}\le C' \|v_\ell\|_{H^1(\Omega_\ell)}
\|v_\ell\|_{L^2(\Omega_\ell)},
\tag{4.3}
$$
with $C'$ independent of $\ell$ (which holds as the domains $\Omega_\ell$
are uniformly Lipschitz).  This allows one to get (4.2) from the arguments
used to prove Theorem 5.1 of [M2], and then (4.1) follows by applying the 
Hodge star operator.

For the convenience of the reader, we describe in more detail how the
Hessian shows up in the proof of this theorem.  The key point is to relate
$\|dv\|^2_{L^2(\Omega_\ell)}+\|\delta v\|^2_{L^2(\Omega_\ell)}$ to the
$H^1$-norm of $v$ for $v$ satisfying $\nu\wedge v=0$, where $\nu=\nu_\ell$.
By density --\,cf., e.g., Proposition 2.5\,-- it suffices to do this analysis 
for $v\in C^2(\Ombar_\ell,\Lambda^*)$. We assume that this is the case
and {\it henceforth drop the subscript} $\ell$. More precisely, let $\wtN$
be any first order differential operator on differential forms with the 
same principal symbol as the Levi-Civita connection $\nabla$.  (The only
reason for not simply taking $\nabla$ is to make the final estimate
depend only on at most the first derivatives of $g$.)  Then $d\delta+\delta d
=-\Delta$ differs from $\wtN^*\wtN$ by a first order operator, so 
$$
(d\delta v,v)_{L^2}+(\delta dv,v)_{L^2}-(\wtN^*\wtN v,v)_{L^2}
=R(v,v),
\tag{4.4}
$$
with $|R(v,v)|\le C\|v\|_{L^2}\|v\|_{H^1}$.  Now, for any first-order
differential operator $P$ on a smooth compact manifold with boundary $X$,
$$
(Pu,v)_{L^2}=(u,P^*v)_{L^2}+\frac{1}{i}\int\limits_{\pa X}
\langle \sigma_P(x,\nu)u,v\rangle\, dS.
\tag{4.5}
$$
As $(1/i)\sigma_d(x,\xi)u=\xi\wedge u,\ (1/i)\sigma_\delta(x,\xi)u=
-\xi\vee u$, and $(1/i)\sigma_\nabla(x,\xi)u=\xi\otimes u$, we deduce that
$$
\aligned
&(d\delta v,v)_{L^2}+(\delta dv,v)_{L^2}-(\wtN^*\wtN v,v)_{L^2} \\
&=\|\delta v\|^2_{L^2}+\|dv\|^2_{L^2}-\|\wtN v\|^2_{L^2}
+\int\limits_{\bW} \Bigl(\langle \nu\wedge \delta v,v\rangle
-\langle \nu\vee dv,v\rangle+\langle\wtN_\nu v,v\rangle\Bigr)\, dS.
\endaligned
\tag{4.6}
$$
Note that $\langle \nu\vee dv,v\rangle=\langle dv,\nu\wedge v\rangle=0$
since $\nu\wedge v=0$, so the middle term in the boundary integral can be
dropped.

To proceed further, extend $\nu$ to a 1-form on a neighborhood of $\bW$,
so $\nu\wedge v$ is a form defined on a neighborhood of $\bW$ as well.
Since $\nu\wedge v$ vanishes on $\bW$, $\nu\wedge v=\rho \tilde{v}$.  
Moreover, 
$$
\delta(f\tilde{v})=f \delta\tilde{v}-df\vee \tilde{v},
\tag{4.7}
$$
so on $\bW$ we have $\langle\delta(\rho\tilde{v}),v\rangle=-\langle
\tilde{v},d\rho\wedge v\rangle=0$, since $d\rho$ is a multiple of $\nu$
and $\nu\wedge v=0$ on $\bW$ by hypothesis.  Thus we can add $\int_{\bW}
\langle \delta(\nu\wedge v),v\rangle\, dS$ to the right side of (4.6),
to obtain
$$
\|\delta v\|^2_{L^2}+\|dv\|^2_{L^2}=\|\wtN v\|^2_{L^2}
-\int\limits_{\bW} \langle \delta(\nu\wedge v)+\nu\wedge \delta v+
\wtN_\nu v,v\rangle\, dS+R(v,v),
\tag{4.8}
$$
with $R(v,v)$ as above.

To examine the integrand in (4.8), consider
$$
P_\nu v=\delta(\nu\wedge v)+\nu\wedge\delta v+\wtN_\nu v.
\tag{4.9}
$$
This is ostensibly a first-order differential operator, but its principal 
symbol satisfies
$$
i \sigma_{P_\nu}(x,\xi)v=\xi\vee(\nu\wedge v)+\nu\wedge(\xi\vee v)
-\langle \nu,\xi\rangle v=0.
\tag{4.10}
$$
Hence $P_\nu$ is actually a zero-order operator.  Moreover, the only term in 
$P_\nu$ that depends on derivatives of $\nu$ is the first one.  Since (4.7)
holds for functions $f$ and forms $\tilde{v}$, if we write 
$\nu=\sum f_i\,dx_i$, then in fact
$$
v\mapsto P_\nu v+\sum_i df_i\vee(dx_i\wedge v)
\tag{4.11}
$$
is not only zero order, but its norm is uniformly bounded as long as the
functions $f_i$ are uniformly bounded on $\bW$, i.e., as long as $\nu$
is uniformly bounded.  Consequently, using
$$
\sum\limits_i\langle df_i\vee(dx_i\wedge v),v\rangle=\sum\limits_i
\langle dx_i\wedge v,df_i\wedge v\rangle=\sum\limits_{i,j}
\frac{\pa^2\rho}{\pa x_i\pa x_j} \langle dx_i\wedge v,dx_j\wedge v\rangle,
\tag{4.12}
$$
we deduce that
$$
\|\delta v\|^2_{L^2}+\|dv\|^2_{L^2}=\|\wtN v\|^2_{L^2}+\int\limits_{\bW}
\sum\limits_{i,j} \frac{\pa^2\rho}{\pa x_i \pa x_j}
\langle dx_i\wedge v,dx_j\wedge v\rangle\, dS+R'(v,v),
\tag{4.13}
$$
with $|R'(v,v)|\le C\|v\|_{L^2} \|v\|_{H^1}$.  For the estimates of the 
theorem, one wants the integral on the right to be positive, modulo terms 
that can be absorbed into $R'$.  This is certainly satisfied for almost
convex domains, and indeed this motivates our definition.  If we add a
large multiple of $\|v\|^2_{L^2}$ to both sides of (4.13), $R'(v,v)$ can be
absorbed by reducing the constants in front of $\|\wtN v\|^2_{L^2}$ and
$\|v\|^2_{L^2}$, giving the desired uniform estimate for $\Omega$ (which, we
recall, stands for $\Omega_\ell$ here).

Parenthetically, we wote that for $(n-1)$-forms $v$,
$$
\langle dx_i\wedge v,dx_j\wedge v\rangle=(dx_i\vee *v)
\overline{(dx_j\vee *v)},
\tag{4.14}
$$
so the quadratic form on normal forms is equivalent to 
$\sum (\pa^2\rho/\pa x_i\pa x_j)\,dx_i\,dx_j$ on the space of vectors 
tangent to $\bW$.

The final ingredient in the proof of Theorem 4.1, as described in the proof
of Theorem 5.1 in [M2], is an approximation argument.  
Consider an approximating sequence $\Omega_\ell\nearrow\Omega$ as
in the proof of Lemma 2.7 and, for 
$u\in \Cal{D}(H^{1/2}_R)=X_R(\Omega)$, let $u_\ell$ be the solution of
$$
(\Delta-1)u_\ell=0,\quad \delta u_\ell=0,\quad 
\nu_\ell\wedge u_\ell=\nu_\ell\wedge u\ \text{ on }\ \pa\Omega_\ell,
\tag{4.15}
$$
satisfying
$$
u_\ell,\ du_\ell\in L^2(\Omega_\ell,\Lambda^*).
\tag{4.16}
$$
Letting $v_\ell=u|_{\Omega_\ell}-u_\ell$, one shows, using (4.13), that 
$v_\ell$ converges to some $v\in H^1(\Omega,\Lambda^*)$ weakly in $H^1$. 
On the other hand, a direct argument using the solution of the auxiliary
problem shows that $u_\ell\rightarrow 0$ weakly in $L^2$; this is quite 
natural since $u_\ell$ solves a homogeneous problem with boundary data
going to $0$ as $\ell\rightarrow\infty$.  Combined, these two show that 
$u=v$, so $u\in H^1(\Omega,\Lambda^*)$.  We refer to [M2] for more
details on this sort of argument.

$\text{}$ \newline
{\smc Remark}. Known results on the Dirichlet and Neumann boundary problems
imply that, when $\Omega$ satisfies the LEBC,
$$
\Cal{D}(H_{A,k}),\ \Cal{D}(H_{R,k})\subset H^2(\Omega,\Lambda^k),\quad
\text{for }\ k=0\ \text{ or }\ n.
$$
However, there is no $\sigma>0$ for which one can say $\Cal{D}(H_{A,1})
\subset H^{1+\sigma}(\Omega,\Lambda^1)$ for all such $\Omega$.  One can see
this by considering the example introduced in \S{1}.  We have
$$
\Cal{D}(H_{A,1})\supset \{df:\,f\in H^2(\Omega),\,\Delta f\in H^1(\Omega),
\pa_\nu f=0\}.
$$
Simple counterexamples show that such $f$ need not belong to 
$H^{2+\sigma}(\Omega)$ for any $\sigma>0$.

$$\text{}$$
{\bf References}

$\text{}$
\roster 
\item"[Ca]" A.~Calder{\'o}n, Cauchy integrals on Lipschitz curves and related
operators, Proc. Nat. Acad. Sci. USA 74 (1977), 1324--1327.
\item"[CMM]" R.~Coifman, A.~McIntosh, and Y.~Meyer, L'int{\'e}grale de Cauchy
definit un op{\'e}rateur born{\'e} sur $L^2$ pour les courbes lipschitziennes,
Annals of Math. 116 (1982), 361--388.
\item"[CD1]" M.~Costabel and M.~Dauge, Un r{\'e}sultat de densit{\'e}
pour les {\'e}quations de Maxwell r{\'e}gularis{\'e}es dans un domaine
lipschitzien, C.~R.~Acad.~Sci.~Paris S{\'e}r. I Math.~327 (1998), 849--854.
\item"[CD2]" M.~Costabel and M.~Dauge, Singularities of electromagnetic fields
in polyhedral domains, Arch. Ration. Mech. Anal. 151 (2000), 221--276.
\item"[F]" K.\,O.~Friedrichs, Differential forms on Riemannian manifolds,
Comm. Pure Appl. Math. 8 (1955), 551--590.
\item"[G]" M.~Gaffney, The harmonic operator for exterior differential
forms, Proc. Nat. Acad. Sci. USA 37 (1951), 48--50.
\item"[MM1]" D.~Mitrea and M.~Mitrea, Finite energy solutions of Maxwell's 
equations and constructive Hodge decompositions on nonsmooth Riemannian
manifolds, J. Funct. Anal. 190 (2002), 339--417.
\item"[MM2]" D.~Mitrea and M.~Mitrea, Poisson problems with half-Dirichlet
boundary conditions for Dirac operators on nonsmooth manifolds,
Mem. on Diff. Eqns. and Math. Phys. 30 (2003), 9--49.
\item"[MMT]" D.~Mitrea, M.~Mitrea, and M.~Taylor, Layer potentials, the
Hodge Laplacian, and global boundary problems in nonsmooth Riemannian 
manifolds, AMS Memoir \#713, Amer. Math. Soc., Providence RI, 2001.
\item"[M1]" M.~Mitrea, Generalized Dirac operators on nonsmooth manifolds
and Maxwell's equations, J. Fourier Anal. Appl. 7 (2001), 207--256.
\item"[M2]" M.~Mitrea, Dirichlet integrals and Gaffney-Friedrichs 
inequalities in convex domains, Forum Math. 13 (2001), 531--567.
\item"[M3]" M.~Mitrea, Sharp Hodge decompositions, Maxwell's equations, and 
vector Poisson problems on non-smooth, three-dimensional Riemannian 
manifolds, Duke Math. J., to appear.
\item"[MT]" M.~Mitrea and M.~Taylor, Boundary layer methods for Lipschitz
domains in Riemannian manifolds, J. Funct. Anal. 163 (1999), 181--251.
\item"[MT2]" M.~Mitrea and M.~Taylor, Potential theory on Lipschitz domains 
in Riemannian manifolds: Sobolev-Besov space results and the Poisson 
problem, J. Funct. Anal. 176 (2000), 1--79.
\item"[Va]" A.~Vasy, Propagation of singularities for the wave equation on
manifolds with corners, Preprint, 2004.
\endroster

\vskip 0.10in
\noindent ------------------------------------------
\vskip 0.08in

\noindent {\tt M.\,Mitrea}: University of Missouri at Columbia, 
Columbia, MO 65211, USA

\vskip 0.06in

\noindent {\tt M.\,Taylor}: UNC at Chapel Hill, Chapel Hill, NC 27599, USA

\vskip 0.06in

\noindent {\tt A.\,Vasy}: Massachusetts Institute of Technology, 
Cambridge, MA 02139, USA

\enddocument